\documentclass[11pt,twoside,a4paper]{article}
\usepackage{amsfonts}
\usepackage{amssymb}
\usepackage{amsmath}

\begin{document}

\newcommand{\pa}{\partial}
\newcommand{\opa}{\overline\pa}
\newcommand{\ol}{\overline }

\numberwithin{equation}{section}

\newcommand\C{\mathbb{C}}  
\newcommand\R{\mathbb{R}}
\newcommand\Z{\mathbb{Z}}
\newcommand\N{\mathbb{N}}
\newcommand\PP{\mathbb{P}}

{\LARGE \centerline{On the Cauchy problem for the $\opa$ operator}}
\vspace{0.8cm}

\centerline{\textsc {Judith Brinkschulte\footnote{Universit\"at Leipzig, Mathematisches Institut, Augustusplatz 10/11, D-04109 Leipzig, Germany. 
E-mail: brinkschulte@math.uni-leipzig.de}
 and C. Denson Hill
\footnote{Department of Mathematics, Stony Brook University, Stony Brook NY 11794, USA. E-mail: dhill@math.sunysb.edu\\
{\bf{Key words:}} Cauchy problem for $\opa$, weakly pseudoconvex hypersurface\\
{\bf{2000 Mathematics Subject Classification:}} 32F10, 32F32, 32W05}}}

\vspace{0.5cm}

\begin{abstract} We present new results concerning the solvability, or lack of thereof, in the Cauchy problem for the $\opa$ operator with initial values assigned on a weakly pseudoconvex hypersurface, and provide illustrative examples.
\end{abstract}

\vspace{0.5cm}

\section{Introduction}

We present here a discussion of an elementary question about the $\opa$ operator. In $\C^n$, or more generally on some $n$-dimensional complex manifold $X$, consider a smooth connected real hypersurface $M$, and a half open - half closed connected domain $D$, which has $M$ as its partial boundary, and is on one side of $M$. Typically, for some connected open set $U$, $D=\lbrace z\in U\mid r(z) \leq 0\rbrace$, $M =\lbrace r(z) =0\rbrace$, where $r: U\longrightarrow \R$ is a $\mathcal{C}^\infty$ function with $dr_{\mid M} \not= 0$. We have in mind the situation where $M$ is {\it not} the full topological boundary of $D$. Any such hypersurface $M$ is non characteristic for the $\opa$ operator, since the Dolbeault complex is an elliptic complex. Also ellipticity and hyperbolicity are not mutually exclusive concepts for complexes of differential operators, as they are for scalar operators: For example the deRham complex is both elliptic and hyperbolic.
 It is reasonable therefore to consider the Cauchy problem for the $\opa$ operator, with assigned "initial values" on $M$. This can be done either at the primitive level of $(p,q)$-forms, or else at the level of cohomology classes. A study of these basic problems was initiated in \cite{AH1} and \cite{AH2}. We recall below the minimal preliminaries needed for this paper, and refer the reader to \cite{AH1} for more details.\\

Let $\mathcal{I}_{(p,q)}(D)$ denote the differential ideal in $\mathcal{C}^\infty_{(p,q)}(D)$ generated by $r$ and $\opa r$; i.e., $u\in \mathcal{I}_{(p,q)}(D)$ means that $u= r\alpha + \opa r\wedge \beta$ where $\alpha\in\mathcal{C}^\infty_{(p,q)}(D)$ and $\beta\in\mathcal{C}^\infty_{(p,q-1)}(D)$. It is intrinsic to the geometry of $D$ and $M$; it does not depend on the choice of the defining function $r$. Smooth tangential $(p,q)$-forms on $M$ are defined as equivalence classes $\lbrack u\rbrack$ in the quotient $\mathcal{Q}_{(p,q)}= \mathcal{C}^\infty_{(p,q)}(D)/ \mathcal{I}_{(p,q)}(D)$, and the tangential Cauchy-Riemann operator $\opa_M : \mathcal{Q}_{(p,q)} \longmapsto \mathcal{Q}_{(p,q+1)}$ is defined by $\opa_M \lbrack u\rbrack = \lbrack \opa u\rbrack$. A form $u\in \mathcal{I}_{(p,q)}(D)$ should be thought of as having zero Cauchy data on $M$. We have the cohomology groups
\begin{eqnarray*}
{} & H^{p,q}(D)&  = \frac{ \mathrm{Ker}\lbrace \opa: \mathcal{C}^\infty_{(p,q)}(D) \rightarrow \mathcal{C}^\infty_{(p,q+1)}(D)\rbrace}{\mathrm{Im}\lbrace \opa: \mathcal{C}^\infty_{(p,q-1)}(D) \rightarrow \mathcal{C}^\infty_{(p,q)}(D)\rbrace},\\
{} & H^{p,q}(D,\mathcal{I}) & =  \frac{ \mathrm{Ker}\lbrace \opa: \mathcal{I}_{(p,q)}(D) \rightarrow \mathcal{I}_{(p,q+1)}(D)\rbrace}{\mathrm{Im}\lbrace \opa: \mathcal{I}_{(p,q-1)}(D) \rightarrow \mathcal{I}_{(p,q)}(D)\rbrace},\\
{} & H^{p,q}(M)&  = \frac{ \mathrm{Ker}\lbrace \opa_M : \mathcal{Q}_{(p,q)}(D) \rightarrow \mathcal{Q}_{(p,q+1)}(D)\rbrace}{\mathrm{Im}\lbrace \opa_M : \mathcal{Q}_{(p,q-1)}(D) \rightarrow \mathcal{Q}_{(p,q)}(D)\rbrace}.
\end{eqnarray*}

Note that our notation is such that all the differential forms involved above are required to be $\mathcal{C}^\infty$ up to the partial boundary $M$ of $D$.\\

Here is the formulation of the Cauchy problem for the $\opa$ operator in terms of differential forms: Given $f\in \mathcal{C}^\infty_{(p,q+1)}(D)$ and $u_0\in\mathcal{C}^\infty_{(p,q)}(M)$, the problem is to find $u\in\mathcal{C}^\infty_{(p,q)}(D)$ such that
\begin{equation} \label{1}
{} \left\{ 
\begin{array}{lll}
 \opa u & = & f\ \mathrm{in}\ D \\
 u_{\mid M} & = & u_0
\end{array}   \right.
\end{equation}

Some compatibility conditions are necessary: we must have $\opa f =0$ on $D$ and $f-\opa \tilde{u}_0\in \mathcal{I}_{(p,q+1)}(D)$ for any $\mathcal{C}^\infty$ extension $\tilde{u}_0$ of $u_0$. The set of solutions to (\ref{1}) is partitioned into equivalence classes by saying that two solutions $u_1$ and $u_2$ are equivalent iff $u_1-u_2$ is cohomologous to zero in $H^{p,q}(D, \mathcal{I})$. Then we have ($\cite[\mathrm{page}\ 351]{AH1}$)
\begin{enumerate}
\item[(i)] The existence of a solution $u\in\mathcal{C}^\infty_{(p,q)}(D)$ to (\ref{1}) for all compatible data $f\in\mathcal{C}^\infty_{(p,q+1)}(D)$ and $u_0\in\mathcal{C}^\infty_{(p,q)}(M)$ is equivalent to $H^{p,q+1}(D,\mathcal{I})=0$.
\item[(ii)] The solution to (\ref{1}) is unique, up to equivalence, if and only if $H^{p,q}(D, \mathcal{I})=0$.
\end{enumerate}

The homogeneous version of the Cauchy problem for the $\opa$ operator is formulated in terms of cohomology classes
as follows: Given a boundary cohomology class $\xi_0\in H^{p,q}(M)$, the problem is to find a cohomology  class $\xi\in H^{p,q}(D)$ such that $\rho (\xi)=\xi_0$, where $\rho$ is the map induced by restriction (dual to smooth extension). Having existence or uniqueness for this version of the Cauchy problem is equivalent to the surjectivity or the injectivity, respectively, of the homomorphism
\begin{equation}  \label{2}
H^{p,q}(D) \overset{\rho}\longrightarrow H^{p,q}(M).
\end{equation}

For this version of the problem it is convenient to also introduce the half open - half closed connected domain $D^+=\lbrace z\in U\mid r(z) \geq 0\rbrace$. Then we have $\cite[\mathrm{page}\ 355]{AH1}$
\begin{enumerate}
\item[(iii)] Either $H^{p,q}(D, \mathcal{I})=0$ or $H^{p,q}(U)=0$ is a sufficient condition for the injectivity (uniqueness) in (\ref{2}).
\item[(iv)] $H^{p,q+1}(D,\mathcal{I})=0$ is a sufficient condition for the surjectivity  (existence) in (\ref{2}), and it is also necessary if $H^{p,q+1}(U)=0$.
\item[(v)] If $H^{p,q}(U)=0$ and $H^{p,q+1}(U)=0$, then for $q>0$, there are isomorphisms
$$H^{p,q}(D^+) \simeq H^{p,q+1}(D,\mathcal{I}),$$
$$H^{p,q}(M) \simeq H^{p,q}(D)\oplus H^{p,q}(D^+).$$
\end{enumerate}

The point of the above discussion is to emphasize that the vanishing, or not, of one of the cohomology groups $H^{p,s}(D,\mathcal{I})$ is an important issue in the Cauchy problem for $\opa$. In particular if one of these cohomology groups is, say, infinite dimensional, then there is a big obstruction to either existence of uniqueness to some Cauchy problem. \\

One might expect that there could be certain finite dimensional obstructions to the existence or uniqueness in the Cauchy problem for the $\opa$ operator. But remarkably, in the situation where $X$ is a Stein manifold, each of the various cohomology groups which enter into the above discussion are either zero or else infinite dimensional $\cite{BHN}$.\\

What happens in the local situation, near a boundary point $x_0\in M$ and for a sufficiently small $D$, was explained in $\cite{AH2}$ in terms of the signature of the Levi form of $M$ at $x_0$, at least under strict assumptions on the number of positive or negative eigenvalues of the Levi form (see Theorems 1,2,...,7 in \cite{AH2}). Since that time a number of authors have obtained results analogous to Theorem 1 in \cite{AH2}, under non-strict assumptions on the Levi-form. But after 35 years there seems to have been little or no progress made in establishing results analogous to Theorem 2 in \cite{AH2}, under non-strict assumptions.\\

In this paper we do obtain results analogous to the above mentioned Theorem 2, under non-strict assumptions, for the special case in which $M$ is weakly pseudoconvex. Our main results (Theorems \ref{first} and \ref{second}) are in fact global in nature. Although we borrow some techniques from \cite{AH2}, it is the considerable technical progress made in \cite{B1} which enables our success.\\

A separate issue is the validity, or nonvalidity, of the Poincar\'e lemma for the tangential Cauchy-Riemann operator $\opa_M$ on $M$. Recently, in a much more general context than $M$ being a hypersurface, some new insight about this question was obtained in \cite{HN}. However this other issue is closely related to the Cauchy problem for the $\opa$ operator; and this is especially so for the case where $M$ is a weakly pseudoconvex hypersurface.\\

\section{The main results}

Throughout this section, we will consider the following set-up:\\

Let $X$ be an $n$-dimensional Stein manifold.
Let $\Omega\subset X$ be an open set, and let $r,g\in\mathcal{C}^\infty(\Omega)$ be two real-valued functions on $\Omega$. We define
\begin{eqnarray*}
M & = & \lbrace z\in\Omega\mid g(z) <0\ \mathrm{and}\ r(z) =0\rbrace,\\
N & = & \lbrace z\in\Omega\mid g(z) =0\ \mathrm{and}\ r(z)\leq 0\rbrace\ \mathrm{and}\\
D & = & \lbrace z\in\Omega\mid g(z)<0\ \mathrm{and}\ r(z)\leq 0\rbrace.
\end{eqnarray*}
We assume that $D$ is connected and relatively compact in $\Omega$. We moreover require that $dr\not= 0 $ on $\ol M$, $dg\not= 0$ on $ N$ and $dr\wedge dg\not= 0$ on $\ol M\cap N$. Our goal is to establish vanishing theorems for the cohomology groups on $D$ with zero Cauchy data on $M$. For this, we need to make the following convexity assumptions on $r$ and $g$:
\begin{eqnarray}  
& i\pa\opa r_{T^{1,0}_x M}\geq 0\ \mathrm{for\ all}\ x\in M
\label{ass1}\\
& i\pa\opa g = 0\ \mathrm{in\ an\ open\ neighborhood}\ V \ \mathrm{of}\  N\ \mathrm{in}\ \Omega
\label{ass2} 
\end{eqnarray}
In particular, the hypothesis on $D$ imply that $D$ is piecewise smooth with a weakly pseudoconvex boundary. With this set-up, we are able to prove the following

\newtheorem{first}{Theorem}[section]
\begin{first}   \label{first} 
We obtain 
$$H^{p,s}(D,\mathcal{I})=0$$
for $p=0,1,\ldots,n$ and $s=0,1,\ldots,n-2$.
\end{first}

{\it Proof.} Let us consider the sets $B(\eta)= D\cap \lbrace -\frac{1}{2} \eta < g < 0\rbrace$ for $0<\eta \leq \eta_0$, where $\eta_0$ is sufficiently small. In particular,  $B(\eta)$ is weakly pseudoconvex.\\

The case $s=0$ is trivial, so let $1\leq s\leq n-2$ and $f\in \mathcal{I}_{(p,s)}(D)$ with $\opa f=0$ be given. Without loss of generality (see Lemma 2.3, page 340, of \cite{AH1}), we may assume that $f$ vanishes to infinite order on $M\cap D$. \\

Let $\chi\in\mathcal{C}^\infty(D)$ be such that
$$ \chi = \left\{ \begin{array}{ccc}
1 & \mathrm{on} & \lbrace g <-\frac{\eta}{3} \rbrace \\
0 & \mathrm{on} & \lbrace g > -\frac{\eta}{4} \rbrace
\end{array} \right. .$$
Then $\chi f \in\mathcal{C}^\infty_{(p,s)}(D)$ and $\opa (\chi f) = \opa\chi\wedge f$ has support in $\lbrace z\in\Omega \mid r(z) \leq 0, -\frac{\eta}{3} \leq g(z) \leq -\frac{\eta}{4}\rbrace \subset \ol B(\eta)$. Using
$\cite[\mathrm{Theorem}\ 4.2]{B1}$, there exists $\varphi\in\mathcal{C}^\infty_{(p,s)}(D)$ such that $\opa\varphi = \opa\chi\wedge\tilde{f}$ and $\mathrm{supp}\varphi\subset \ol B(\eta)$. We set $f_1 =\chi f -\varphi$. Then clearly $\opa f_1 =0$ and $f_1 =f$ on $D\cap\lbrace g< -\frac{\eta}{2}\rbrace$. Moreover, $\mathrm{supp} f_1 \subset \ol D$.  Hence we may again apply $\cite[\mathrm{Theorem}\ 4.2]{B1}$ in order to obtain a solution $u\in\mathcal{I}_{(p,s-1)}(D)$ to the equation $\opa u = f_1$.\\

We conclude that given $f\in\mathcal{I}_{(p,s)}(D)$ with $\opa f=0$, $1\leq s\leq n-2$, and given $\eta >0$, we can find $u_\eta \in\mathcal{I}_{(p,s-1)}(D_\eta)$ such that $\opa u_\eta =f$ on $D_\eta$, where $D_\eta = D\cap \lbrace g < -\eta\rbrace$.\\
 
Now consider a strictly decreasing sequence $(\eta_\nu)_{\nu\in\N}$ with $\eta_\nu \underset{\nu\rightarrow\infty} \rightarrow 0$. Then, for each $\nu\in\N$, we can find $u_\nu\in \mathcal{I}_{(p,s-1)}(D_{\eta_\nu})$ such that $\opa u_\nu = f$ on $D_{\eta_\nu}$. Therefore $\opa (u_\nu - u_{\nu-1})=0$ on $D_{\eta_{\nu-1}}$. If $s=1$, then $u_2=u_1$ on $D_{\eta_1}$, $u_3=u_2$ on $D_{\eta_2}, \ldots.$ Hence, by setting $u=u_\nu$ on $D_{\eta_\nu}$, we define an element $u\in \mathcal{I}_{(p,s-1)}(D)$ with $\opa u =f$.\\

If $s>1$, then we construct a sequence $(u^\prime_\nu)_{\nu\in\N}$, $u_\nu^\prime\in\mathcal{I}_{(p,s-1)}(D_{\eta_\nu})$ satisfying $\opa u^\prime_\nu = f$ on $D_{\eta_\nu}$, $u^\prime_\nu = u^\prime_{\nu-1}$ on $D_{\eta_{\nu-3}}$. Indeed, suppose $u^\prime_1,\ldots, u^\prime_\nu$ have been constructed. Then $\opa (u_{\nu+1}-u_\nu^\prime)=0$ on $D_{\eta_\nu}$. Since the domain $D_{\eta_\nu}$ satisfies the same hypotheses as $D$, there exists $\sigma\in \mathcal{I}_{(p,s-2)}(D_{\eta_{\nu-1}})$ such that $\opa\sigma = u_{\nu+1} - u^\prime_\nu$ on $D_{\eta_{\nu-1}}$. Let $\tau$ be a smooth function on $D$ with
$$\tau = \left\{ \begin{array}{cc}
1 & \mathrm{on}\ D_{\eta_{\nu-2}}\\
0 & \mathrm{outside}\ D_{\eta_{\nu-1}}
\end{array} \right. $$
Set $u^\prime_{\nu+1}= u_{\nu+1}-\opa(\tau\sigma)$. Then $u^\prime_{\nu+1}$ has the required properties. Now setting $u=u^\prime_\nu$ on $D_{\eta_{\nu-2}}$, we get a well-defined element $u\in \mathcal{I}_{(p,s-1)}(D)$ for which $\opa u =f$. \hfill$\square$\\

Now let us replace the assumptions (\ref{ass1}) and (\ref{ass2}) by the following more restrictive assumptions on $r$ and $g$

\begin{eqnarray}  
i\pa\opa r \geq 0\ &\mathrm{in\ an\ open\ neighborhood}\ U\ \mathrm{of}\ \ol M\ \mathrm{in}\ \Omega
\label{ass3}\\
i\pa\opa g = 0\ &\mathrm{in\ an\ open\ neighborhood}\ V \ \mathrm{of}\ \ol D\ \mathrm{in}\ \Omega
\label{ass4} \\
i\pa\opa r >0 \ &\mathrm{in\ an\ open\ neighborhood}\ W \ \mathrm{of}\ \ol M\cap N\ \mathrm{in}\ \Omega
\label{ass5}
\end{eqnarray}

We then obtain the following result.

\newtheorem{second}[first]{Theorem}
\begin{second} \label{second}
We obtain
$$H^{p,s}(D,\mathcal{I})=0$$
for $p=0,1,\ldots,n$ and $s=0,1,\ldots,n-1$.
\end{second}

Before we start the proof of the theorem, let us consider the domains
\begin{eqnarray*}
A(\varepsilon) & = & \lbrace z\in\Omega\mid r(z) <\varepsilon,\ g(z) <0\rbrace,\\
B(\varepsilon,\eta) & = & \lbrace z\in\Omega\mid r(z) <\varepsilon,\ -\frac{1}{2}\eta < g(z) < 0\rbrace
\end{eqnarray*}
for $0<\varepsilon\leq \varepsilon_0$, $0<\eta\leq\eta_0$, where $\varepsilon_0 ,\eta_0$ are sufficiently small. In particular, both $A(\varepsilon)$ and $B(\varepsilon, \eta)$ are weakly pseudoconvex domains.

\newtheorem{runge}[first]{Lemma}
\begin{runge}  \label{runge} 
$(A(\varepsilon), B(\varepsilon,\eta))$ is a Runge pair in all degrees $s\leq n$, i.e. the natural map
$$H^{p,s}_c(B(\varepsilon,\eta))\longrightarrow H^{p,s}_c((A(\varepsilon))$$
is injective for all $0\leq p,s\leq n$. Here $H^{p,s}_c(B(\varepsilon,\eta))$ and $H^{p,s}_c((A(\varepsilon))$ denote the $\opa$-cohomology groups for smooth forms with compact support in $B(\varepsilon,\eta)$ resp. $A(\varepsilon)$.
\end{runge}

{\it Proof of the Lemma.} By a criterion for Runge pairs proved in \cite[p. 122]{AV}, for each compact $K\subset B(\varepsilon,\eta)$, it suffices to construct a smooth, strictly plurisubharmonic exhaustion function $\Phi$ on $A(\varepsilon)$ such that
$$K\subset \lbrace z\in A(\varepsilon)\mid \Phi(z) \leq \sup_K\Phi\rbrace \subset B(\varepsilon,\eta).$$
Consider the three plurisubharmonic functions
\begin{eqnarray*}
\varphi_1 & = & r-\varepsilon,\\
\varphi_2 & = & g,\\
\varphi_3 & = & -g-\frac{1}{2}\eta
\end{eqnarray*}
Note that it is no loss of generality to assume that $r$ is plurisubharmonic on an open neighborhood of $\ol D$ in $\Omega$ (Indeed, replacing $r$ by $\max (-\tau,r)$ for some small $\tau >0$, one obtains a continuous plurisubharmonic function in a neighborhood of $\ol D$ in $\Omega$ that is still a defining function for $M$). Then $\varphi = \max(\varphi_1,\varphi_2,\varphi_3)$ is a continuous plurisubharmonic function on $A(\varepsilon)$, and $B(\varepsilon,\eta) = \lbrace z\in A(\varepsilon)\mid \varphi(z) <0\rbrace$. Let $K\subset B(\varepsilon,\eta)$ be compact. It is no loss of generality to assume $K=\lbrace z\in A(\varepsilon)\mid \varphi(z) < -\delta\rbrace$ for some $\delta > 0$. But since $\varphi$ is plurisubharmonic on $A(\varepsilon)$, it then follows by definition of the plurisubharmonic hulls that $K = \hat{K} {}_{B(\varepsilon,\eta)}= \hat{K} {}_{A(\varepsilon)}$. So there exists a smooth plurisubharmonic exhaustion function $\Phi$ on $A(\varepsilon)$ with $\Phi < 0$ on $K$ and $\Phi > 0$ outside $B(\varepsilon,\eta)$. \hfill$\square$\\

{\it Proof of the theorem.} Let $1\leq s\leq n-1$ and $f\in \mathcal{I}_{(p,s)}(D)$ with $\opa f=0$ be given. Without loss of generality (see Lemma 2.3, page 340, of \cite{AH1}), we may assume that $f$ vanishes to infinite order on $M\cap D$. We extend $f$ to $\tilde{f}$, defined on $A(\varepsilon)$, by defining $\tilde{f}$ to be zero outside of $D$.\\

Let $\chi\in\mathcal{C}^\infty(A(\varepsilon))$ be such that
$$ \chi = \left\{ \begin{array}{ccc}
1 & \mathrm{on} & \lbrace g <-\frac{\eta}{3} \rbrace \\
0 & \mathrm{on} & \lbrace g > -\frac{\eta}{4} \rbrace
\end{array} \right.$$
Then $\chi\tilde{f} \in\mathcal{C}^\infty_{0(p,s)}(A(\varepsilon))$ and $\opa (\chi\tilde{f}) = \opa\chi\wedge\tilde{f}$ has support in $\lbrace z\in\Omega \mid r(z) \leq 0, -\frac{\eta}{3} \leq g(z) \leq -\frac{\eta}{4}\rbrace \subset\subset B(\varepsilon,\eta)$. Using the lemma above, there exists $\varphi\in\mathcal{C}^\infty_{0(p,s)}(B(\varepsilon,\eta))$ such that $\opa\varphi = \opa\chi\wedge\tilde{f}$. We set $f_1 =\chi\tilde{f} -\varphi$. Then clearly $\opa f_1 =0$ and $f_1 =f$ on $D\cap\lbrace g< -\frac{\eta}{2}\rbrace$.\\

Now let $\rho$ be a smooth function on $\ol\Omega$ with $0\leq\rho\leq 1$ such that
$$\rho = \left\{ \begin{array}{ccc}
1 & \mathrm{on} & \lbrace g\geq - \frac{\eta}{2} \rbrace\\
0 & \mathrm{on} & \lbrace g \leq -\eta\rbrace,
\end{array} \right.$$
and let
$$ D^\prime (\varepsilon,\eta)= \lbrace r\leq \varepsilon\rho\rbrace\cap \lbrace g <0\rbrace.$$
Then $\mathrm{supp} f_1\subset \ol{D^\prime(\varepsilon,\eta)}$. Moreover, using (\ref{ass3}), the bumped domain $D^\prime(\varepsilon,\eta)$ is piecewise smooth with a weakly pseudoconvex boundary if $\varepsilon$ and $\eta$ are sufficiently small. Hence we may apply the results of $\cite[\mathrm{Theorem}\ 4.2]{B1}$ to the domain $D^\prime(\varepsilon,\eta)$ in order to obtain a solution $u\in\mathcal{I}_{(p,s-1)}(D^\prime(\varepsilon,\eta))$ to the equation $\opa u = f_1$.\\

As before, we conclude that given $f\in\mathcal{I}_{(p,s)}(D)$ with $\opa f=0$, $1\leq s\leq n-1$, and given $\eta >0$, we can find $u_\eta \in\mathcal{I}_{(p,s-1)}(D_\eta)$ such that $\opa u_\eta =f$ on $D_\eta$, where $D_\eta = D\cap \lbrace g < -\eta\rbrace$. The rest of the proof is the same as the proof of Theorem \ref{first}.\hfill$\square$\\

As explained in the introduction, Theorems \ref{first} and \ref{second} have applications to the solvability of the tangential Cauchy-Riemann equation on $M$ without shrinking the domain. For the case of $M$ being a pseudo-convex $CR$ hypersurface of finite type and $N$  being a flat hypersurface, it was proved in \cite{S} that $H^{p,q}(M)=0$ for $1\leq q < n-2$. In \cite{LT} and \cite{FLT}, one can find cohomological and geometrical characterizations of the open subsets of a strictly pseudoconvex boundary in a Stein manifold on which one can solve the tangential Cauchy-Riemann equation in all bidegrees without shrinking the domain. Our main results imply the following

\newtheorem{tang}[first]{Corollary}
\begin{tang}  \label{tang} 
Assume that $Y=\lbrace z\in\Omega\mid g(z) > 0\rbrace$ is Stein and that (\ref{ass3}) is satisfied.
\begin{enumerate}
\item If (\ref{ass2}) holds, then we have
$$H^{p,q}(M)=0$$
for $0\leq p\leq n, 1\leq q < n-2$.
\item If (\ref{ass4}) and (\ref{ass5}) hold, then we have
$$H^{p,q}(M)=0$$
for $0\leq p\leq n, 1\leq q <n-1$.
\end{enumerate}
\end{tang}

{\it Proof.} The assumption that $Y$ is Stein implies that
$$H^{p,q}(M) \simeq H^{p,q}(D) \oplus H^{p,q+1}(D,\mathcal{I})$$
(see (v) in the introduction, with $U=Y$). If (\ref{ass3}) is satisfied, we may moreover apply Dufresnoy's results \cite{D} on the solvability of $\opa$ with regularity up to the boundary on weakly pseudoconvex domains. Applying a standard Mittag-Leffler-type procedure, one can conclude that for the half open - half closed domain $D$ one has $H^{p,q}(D)=0$ for $1\leq q\leq n$. We conclude by evoking Theorem \ref{first} and \ref{second}. \hfill$\square$\\

{\it Example.} Here is a simple but illustrative example. Consider the unit sphere
$$S^{2n-1}: \ \vert z_1\vert^1 + \vert z_2\vert^2 + \ldots + \vert z_n\vert^2 =1,$$
which is the boundary of the closed unit ball $\ol B$ in $\C^n$, $n\geq 3$. Let $\tilde{D}$ be the set $\ol B\cap\lbrace \eta < x_1 \leq \eta +\varepsilon\rbrace$, where $0<\eta<\eta +\varepsilon < 1$. It has the partial boundary $\lbrack S^{2n-1}\cap \lbrace \eta < x_1\leq \eta + \varepsilon\rbrace\rbrack \cup \lbrack \ol B \cap \lbrace x_1 \equiv \eta + \varepsilon\rbrace\rbrack$. We make a small $\mathcal{C}^\infty$ "rounding off of the corners" of the non-smooth part of the partial boundary of $\tilde{D}$, and let $D$ denote the resulting half open - half closed domain, and $M$ its smooth partial boundary. This can be done in such a way that (\ref{ass3}) is satisfied. Then from Theorem \ref{second} and Corollary \ref{tang} we obtain that, for $p=0,1,\ldots,n$,
\begin{equation}  \label{6}
\begin{array}{lll}
H^{p,s}(D,\mathcal{I})=0, & & \mathrm{when}\ s =0,1,\ldots, n-1,\\
H^{p,q}(M)=0, & & \mathrm{when}\ q = 1,2,\ldots, n-2.
\end{array}
\end{equation}

Note that here $\varepsilon >0$ is arbitrarily small. Hence $M$ {\it is Levi-flat except for an annular ring having arbitrarily small measure}, where $M$ is strictly pseudoconvex.\\

Thus, in this example, the only homogeneous Cauchy problems (\ref{2}) which are of any interest are those with $q=0$ and $q=n-1$. They both have infinite dimensional spaces of Cauchy data. For $q=0$, there is the well-known isomorphism $H^{p,0}(D)\simeq H^{p,0}(M)$. But for $q=n-1$, the well-posed Cauchy problem occurs on the other side: $H^{p,n-1}(D^+)\simeq H^{p,n-1}(M)$, where $D^+=\lbrace \eta < x_1\rbrace \setminus\overset{\circ}D$. Here we use again (v) and the fact that $Y=\lbrace \eta < x_1\rbrace$ is Stein. Note that the cohomology groups which are missing in (\ref{6}), namely $H^{p,n}(D,\mathcal{I}), H^{p,0}(M)$ and $H^{p,n-1}(M)$ are all {\it infinite dimensional} (see \cite{AH2}).\\

\section{How not to slice eggs}

Denote the coordinates in $\C^n$ by $(z,w)$ with $z=(z_1,z_2)$ and $w=(w_1,w_2,\ldots, w_{n-2})$, and consider
\begin{equation}  \label{a}
\pa\Omega : \vert z_1\vert^2 + \vert z_2\vert^2 + \vert w_1\vert^{m_1} + \ldots + \vert w_{n-2}\vert^{m_{n-2}}=1,
\end{equation}
where $m_1, m_2,\ldots, m_{n-2}$ are even integers, all $\geq 4$. Then $\pa\Omega$ is the weakly pseudoconvex boundary of a generalized closed convex egg $\ol\Omega$ in $\C^n$. For $j=0,1,\ldots, n-2$ let $\Sigma_{n-1-j}$ be the set of points on $\pa\Omega$ at which {\it exactly} $j$ components of $w$ are zero. Then $\pa\Omega=\cup_{k=1}^{n-1}\Sigma_k$. At each point $x_0\in\Sigma_k$ the complex hessian of $r=\vert z_1\vert^2 + \vert z_2\vert^2 + \vert w_1\vert^{m_1} + \ldots + \vert w_{n-2}\vert^{m_{n-2}}-1$ has $k+1$ positive and $n-k-1$ zero eigenvalues. Hence the Levi form of $\pa\Omega$ at $x_o$ has $k$ positive and $n-k-1$ zero eigenvalues. Note that $\Sigma_k$ has real codimension $2(n-1)-2k$ in $M$, and this codimension is equal to $2j$ when $k=n-1-j$. Thus the real codimension in $M$ of the locus of degeneracy of the Levi form is equal to two times the number of components of $w$ which have been set equal to zero.\\

Now fix a $k$ $(1\leq k < n-1)$, take some point $x_0\in \Sigma_k$, and slice $\ol\Omega$ by a hyperplane which is parallel to $T_{x_0}\pa\Omega$, so as to cut out a sufficiently small half open - half closed convex domain $D$, with a partial boundary $M$. Then the Levi form of $M$ has at least $k+1$ positive eigenvalues at each point, except along the locus $\Sigma_M= M\cap\Sigma_k$, where there are only $k$ positive eigenvalues. This $\Sigma_M$ is {\it thin}, having real codimension $2n-2-2k$ in $M$, but it is {\it not compact} in $M$. Hence we cannot apply part 2 of Corollary \ref{tang}, but we may apply part 1. The conclusion is that for $0\leq p\leq n$ we have
\begin{equation}  \label{b}
H^{p,q}(M)=0,\ \mathrm{when}\ 1\leq q < n-2.
\end{equation}
Note that the result is independent of the choice of $k$. One could argue that this is the {\it usual way} to slice an egg. It enables us, in bidegree $(p,q)$, to solve the tangential Cauchy-Riemann equations on $M$ without shrinking.\\

Next we slice the egg in a {\it different way}: with the same choice of the point $x_0$ as before, fix any Riemannian metric $g$ on $\pa\Omega$ and denote by $B(x_0,r)$ the open ball on $\pa\Omega$, centered at $x_0$, of radius $r$. For example, we may use the standard metric on $\pa\Omega$ which is induced by the euclidean metric in the ambient $\C^n$. It was shown in \cite{HN} (see Theorem 7.2, and page 218) that if $r>0$ is taken to be sufficiently small, then
\begin{equation}  \label{c}
\dim H^{p,k}(B(x_0,r))=\infty, \ \mathrm{for\ all} \ 0\leq p\leq n.
\end{equation}

The contrast between (\ref{b}) and (\ref{c}) came, at first, as a bit of a surprise to the authors. However we are able to explain it by means of some subtle but elementary geometry: let us write the tangential Cauchy-Riemann equations on $\pa\Omega$ as
\begin{equation} \label{d}
\opa_M u=f,
\end{equation}
\begin{equation}  \label{e}
\opa_M f=0.
\end{equation}
Let $f$ be a tangential $(p,k)$ form which is $\mathcal{C}^\infty$ and satisfies (\ref{e}) on $B(x_0,r)$. We seek a tangential $(p,k-1)$ form $u$ that is $\mathcal{C}^\infty$ and solves (\ref{d}) on $B(x_0,r^\prime)$, for some $r^\prime$, with $0< r^\prime\leq r$ (shrinking is now allowed). What was shown in \cite{HN} is that there exist constants $r_0 > 0$ and $C>0$ such that, for $0< r^\prime \leq r \leq r_0$, if $r^\prime > C r^{3/2}$, then there is an infinite dimensional space of such $f$'s to which there does not correspond any such solution $u$. The infinite dimensionality expressed in (\ref{c}) then follows.\\

In order to most simply explain why (\ref{b}) and (\ref{c}) can both be true, let us work with the most dangerous case in which all $m_j=4$ in (\ref{a}). When we slice the egg $\Omega$ in the usual way, the $M$ which is cut out looks like a curved "elliptical" surface: there is a longest geodesic on $\pa\Omega$ from $x_0$ to the boundary of $M$, and there is a shortest geodesic on $\pa\Omega$ from $x_0$ to the boundary of $M$. Let $r$ be the length of the longest geodesic, and $r^\prime$ be the length of the shortest geodesic. Denote by $h>0$ the perpendicular distance between the hyperplane $T_{x_0}\pa\Omega$ and the parallel hyperplane used to cut out the slice $D$ of the egg $\Omega$. Thus $r$ and $r^\prime$ are functions of $h$. As $h\rightarrow 0$ we need to obtain the asymptotic behavior of the ratio of two arc-length integrals. For the shorter one it is sufficient to consider $x^2 + y^2 =1$ and hence
$$ r^\prime (h)= \int_{1-h}^1 \frac{dx}{\sqrt{1-x^2}}.$$
For the longer one it suffices to consider $x^2 + y^4=1$ and therefore
$$r(h)= \int_{1-h}^1 \frac{\lbrack 4(1-x^2)^{3/2}+ x^2\rbrack^{1/2}}{2\lbrack 1-x^2\rbrack^{3/4}} dx.$$
A computation shows that
$$\frac{r^\prime(h)}{\lbrack r(h)\rbrack^{3/2}} = \mathrm{const}\ h^{1/8} + O(h^{9/8}),\ \mathrm{as} \ h\rightarrow 0.$$
This means that for $h$ taken sufficiently small, we cannot maintain the crucial inequality $\frac{r^\prime}{r^{3/2}} > C$ from \cite{HN}. Thus (\ref{b}), obtained by slicing the egg in the usual way, does not contradict (\ref{c}) being valid, which was obtained actually not by slicing, but by taking small balls in some Riemannian metric on the boundary of the egg.\\


\begin{thebibliography}
\footnotesize


\bibitem[AH1]{AH1} \textsc{A. Andreotti, C.D. Hill:} \emph{E.E. Levi Convexity and the Hans Lewy problem. Part I.} Ann. Sc. Norm. Super. Pisa {\bf 26}, 325--363 (1972).

\bibitem[AH2]{AH2} \textsc{A. Andreotti, C.D. Hill:} \emph{E.E. Levi Convexity and the Hans Lewy problem. Part II.} Ann. Sc. Norm. Super. Pisa {\bf 28}, 747--806 (1972).

\bibitem[AV]{AV} \textsc{A. Andreotti, E. Vesentini:} \emph{Carleman estimates for the Laplace-Beltrami equation in complex manifolds.} Publ. Math. I.H.E.S. {\bf 25}, 81--130 (1965).

\bibitem[B]{B1} \textsc{J. Brinkschulte:} \emph{The $\opa$-problem with support conditions on some weakly pseudoconvex domains.} Ark. Mat. {\bf 42}, 259--282 (2004).

\bibitem[BHN]{BHN} \textsc{J. Brinkschulte, C.D. Hill, M. Nacinovich:} \emph{Obstructions to generic embeddings.} Ann. Inst. Fourier {\bf 52}, 1785--1792 (2002).

\bibitem[D]{D} \textsc{A. Dufresnoy:} \emph{Sur l'op\'erateur $\opa$ et les fonctions diff\'erentiables au sens de Whitney.} Ann. Inst. Fourier {\bf 29}, 229--238 (1979).

\bibitem[FLT]{FLT} \textsc{F. Forstneric, Ch. Laurent-Thi\`ebaut:} \emph{Stein compacts in Levi-flat hypersurfaces.} Pr\'epubl. Inst. Fourier {\bf 659} (2004).

\bibitem[HN]{HN} \textsc{C.D. Hill, M. Nacinovich:} \emph{On the failure of the Poincar\'e lemma for $\opa_M$ II}, Math. Ann. {\bf 335}, 193--219 (2006).

\bibitem[LT]{LT} \textsc{Ch. Laurent-Thi\'ebaut:} \emph{Sur l'\'equation de Cauchy-Riemann tangentielle dans une calotte strictement pseudoconvexe.} Int. J. Math. {\bf 16}, 1063--1079 (2005).

\bibitem[N]{N} \textsc{M. Nacinovich:} \emph{On boundary Hilbert differential complexes.} Annales Polonici Mathematici {\bf XLVI} 213--235 (1985).

\bibitem[S1]{S} \textsc{M.-C. Shaw:} \emph{Local existence theorems with estimates for $\opa_b$ on weakly pseudo-convex $CR$ manifolds.} Math. Ann. {\bf 294}, 677--700 (1992).

\bibitem[S2]{S2} \textsc{M.-C. Shaw:} \emph{Semi-Global Existence Theorems of $\opa_b$ for (0,n-2) forms on Pseudo-convex boundaries in $C^n$.} Ast\'erisque, Soci\'et\'e Math\'ematique de France, Colloque d'Analyse complexe et g\'eometrie, Marseille (1993), 227--240.


\end{thebibliography}
\end{document}